\documentclass{siamltex}   % Specifies the document style.

\def\Box{{\hfill\hbox{\enspace${\square}$}} \smallskip}
\def\sqr#1#2{{\vcenter{\vbox{\hrule height .#2pt
                             \hbox{\vrule width .#2pt height#1pt \kern#1pt
                                   \vrule width .#2pt}
                             \hrule height .#2pt}}}}
\def\square{\mathchoice\sqr54\sqr54\sqr{4.1}3\sqr{3.5}3}

\newcommand{\be}[1]{\begin{equation}\label{#1}}
\newcommand{\ee}{\end{equation}}
\newcommand{\bl}[1]{\begin{lemma}\label{#1}}
\newcommand{\bp}[1]{\begin{proposition}\label{#1}}
\newcommand{\br}[1]{\begin{remark}\label{#1}}
\title{Rigorous WKB for finite  order Linear
Recurrence Relations with Smooth Coefficients
\thanks{Work supported, in part, by AFOSR Grant AF-0115}}
\author{Ovidiu Costin$^{^{\dag}}$ and  Rodica Costin%
\thanks {Mathematics Department, Hill Center
Rutgers University, New Brunswick, NJ 08903;
e-mail: COSTIN@MAXWELL.RUTGERS.EDU}}
\begin{document}          
\def\z{\noindent}
\def\ep{\epsilon}
\def\epi{\epsilon ^{-1}}
\def\Nsm{\hbox{\small I\hskip -2pt N}}
\def\rdb{\hbox{ I\hskip -2pt R}}
\def\cdb{\hbox{\it l\hskip -5.5pt C\/}}
\def\ndd{\hbox{\it I\hskip -2pt N}}
\def\zdd{\hbox{\sf Z\hskip -4pt Z}}
\def\e{\epsilon}
\def\ei{\epsilon^{-1}}
\def\a{\alpha}
\def\l{\lambda}
\def\tl{\tilde\lambda}
\def\ct{{\mbox{const}}}
\def\d{\delta}
\def\cf{{\cal F}}
\maketitle
\hskip -0.5truecm\begin{abstract} We study the $\epsilon \rightarrow 0$ behavior of 
recurrence relations of the type
$\sum_{j=0}^l a_j(k\epsilon,\epsilon)y_{k+j}=0,$ 
$k\in \zdd$ ($l$ fixed). The $a_j$ are $C^{\infty}$
functions in each variable on $I\times [0,\e_0]$ for a bounded
interval $I$ and $\e_0>0$. Under certain regularity assumptions we find 
the asymptotic behavior of the solutions of such recurrences.
In typical cases 
there exists a fundamental set of solutions in the form
$\{\exp(\epi F_m(k\epsilon,\epsilon))\}_{m=1\ldots l}$ where
the functions $F_m$ are $C^{\infty}$ in each variable on
the same domain as the $a_j$, showing in particular that
the formal  perturbation-series solutions are asymptotic to true
solutions of these recurrences. 
 Some applications are also briefly discussed.
\end{abstract}
%%\\
%%-----------------------------------------\\
%%$^*${\small{Work supported, in part, by AFOSR Grant AF-0115}}$\\$
%%$^1${\small{e-mail COSTIN@MAXWELL.RUTGERS.EDU}}$\\$
%%$^2${\small{e-mail RCOSTIN@MATH.RUTGERS.EDU}}$\\$
{\small\bf\hskip -0.2 truecm AMS(MOS) subject classifications:} 41A60, 65M06, 65M12
\section{Introduction}In the present paper we study the asymptotic
behavior to all orders in $\e\rightarrow 0$ 
of the solutions of one-dimensional recurrence relations of the
form
\be{mainrec}
\sum_{j=0}^l a_j(k\epsilon ,\epsilon )y_{k+j}=0 
\ee
which we may interpret as follows: for each
fixed $k,\ y_{k+l}$ is determined from its predecessors
 $y_{k}\ldots y_{k+l-1}$ (this is assumed possible ---  see  
condition (\ref{asnonz}) below).

Under some further
regularity assumptions we prove that 
the general solution of the recurrence can be piecewise
represented as a sum

\be{expWKB}
y_{k,\e}=\sum_{m=1}^lC_{m,p} \exp[\epi F_m(k\epsilon ,\epsilon )]
\ee
where the functions $F_m$ are everywhere smooth with the
exception of a small neighborhood of the points
where two characteristic roots (\ref{charpol}) cross and where
the representation is different (Proposition 2.2) . In particular,
a fundamental system of solutions can be chosen such that
each of them has, for small--$\e$ and between crossings, a WKB--like expansion:

\be{WKB}
y_k\sim\exp\left(\e\phi(k\e)\right)\left(A_0(k\e)+\e
A_1(k\e)+...\right)
\ee
where $\phi$ is the root of the eikonal equation
$$\exp(\phi'(x))=\lambda(x)$$

and $\lambda$ is one of the $l$ roots of the characteristic
equation  (\ref{charpol}) and the successive amplitudes $A_i$
can be determined by perturbation--expansion.  For technical reasons, we prefer the less familiar
notation $(\ref{expWKB})$. 
It is essential for our arguments that a continuous branch of $\ln\,A_0$
can be chosen.

One of the  applications of the rigorous WKB--approach
for discrete schemes is in  determining the spectrum of
large matrices with slowly--varying entries. Such problems
appear for instance in the continuum limit of the Toda lattice.
This system , described by the Hamiltonian 
$H=\frac{1}{2}\sum_{k=1}^{n}p_k^2+\sum_{k=1}^{n-1}\exp(x_k-x_{k-1})$
is completely integrable; this can be expressed in
terms of the constancy of the spectrum of the matrix 

\be{todam}
\left(
\begin{array}{ccccc}
a_1 & b_1 & 0 & ... & 0 \\
b_1 & a_2 & b_2 & ...&  0 \\
0 & b_2 & a_3 & ...& 0 \\
... & ... & ... & ... & ... \\
... & ... & ... & ... & ... \end{array}
\right)
\ee

The spectral problem in case the coefficients
are of the form $a_k=a(k\e),\ b_k=b(k\e)$ with
$a$ and $b$ smooth and satisfying some regularity
conditions leads to a recurrence
of the type (\ref{mainrec}) that is solved asymptotically 
by the methods described below [1].
\vskip 0.5 cm
In  (\ref{mainrec}) the number $l$ of steps  of the recurrence
is fixed, $\ep > 0$ is
a small parameter and $k \in
\zdd$ is such that $k\ep \in { I}$ where ${ I}\subset
\rdb$ is a compact interval. Some initial or boundary conditions
are assumed.

The coefficients
$a_0(x,\e)\ldots a_l(x,\e): I\times[0,\e_0]\rightarrow
\cdb$ are assumed $C^{\infty}$ in  $x$ and in $\e$ in some
domain $I\times[0,\e_0]$. We also require the existence of a uniformly
asymptotic series for $a_j$: for any $t\in\ndd$,

\be{asona}
\mid a_j(x,\e)-\sum_{s=0}^t a_{j,s}(x)\e^s \mid <M_{j,s}\e^{t+1}
\ee
where the functions $a_{j,s}$ are 
$C^{\infty}$ in $x$ 
(for instance $a_j \in C^\infty( I $ $\times[0,\e_0]$).

We are also imposing the nonsingularity condition 
\be{asnonz}
\inf_{x \in I}\{|a_0(x,0)|\}>0
{\mbox{ and}\ \ }\inf_{x \in I}\{| a_l(x,0)|\}>0
\ee

We begin by giving some simple examples and deriving heuristically
their small-$\ep$ behavior. The contents of the paper
will subsequently make the given solutions rigorous.

{\bf a.} Consider the one-step recurrence relation
\be{reuler}
y_{k+1}=e^{f(k\epsilon )}y_k
\ee
where  $y_0=1$, $0\leq k\leq \epsilon^{-1}$ and $f$ is a $C^\infty$ 
function on  $[0,1]$.
It has the explicit solution
\be{expeu}
y_{k}=\exp((\sum_{j=0}^{k-1} f(j\epsilon ))
\ee

When $f=f_0$ is constant, $y_{k}=\exp(\epi {f_0}k\ep)$ and this should
also be the  order of magnitude of $y_{k }$ for a general smooth f
when $\ep$ is very small;
we then try a formal asymptotic solution
\be{formalseries}
y_k\sim~\exp(\epsilon^{-1}\Phi_0(k\epsilon)+\Phi_1(k\epsilon)+\ep\
\Phi_2(k\ep)+\ldots )
\ee
 in (\ref{reuler}) . To be consistent with (\ref{reuler}) we must have
\be{consistc}
\exp[\sum_{m=0}^{\infty}\ep^{m-1}\Phi_m(k\ep+\ep)]\sim
\exp[f(k\ep)+\sum_{m=0}^{\infty}\ep^{m-1}\Phi_m(k\ep)],\  \ep \rightarrow 0
\ee
Expanding the exponent in the LHS of (\ref{consistc}) in a Taylor series
around $k\ep$ and then 
identifying the corresponding powers of $\ep$ in  (\ref{consistc}) we
obtain 
\be{indueu}
\Phi'_0=f,\ \ \ \ \Phi'_1=-\frac{1}{2}f' \ \ldots \ \ 
\Phi'_j=C_j\ f^{(j)},\ldots 
\ee
where the constants $C_j$ are seen to be $f$-independent and can thus
be determined  by choosing a particular  $f$ 
for which the sum in (\ref{expeu}) can be explicitly evaluated.
Taking, e.g.,  $f(x)=\exp( x)$  this  sum is 
%%%
%%% 
\be{defbern}
(\exp (x)-1)(\exp (\ep)-1)^{-1} \sim \ep^{-1}(\exp (x)-1) \sum_{j=0}^{\infty}
\frac{B_{2j}}{(2j)!}\ep ^{2j},\ \ \ep \rightarrow 0
\ee
where $B_n$ are the Bernoulli numbers. 
It follows that $C_{j}=B_{j}$ ($B_{2j+1}=0$) and thus, for any
smooth $f$
\be{eumclr}
\sum_{j=0}^{k-1} f(j\epsilon )\sim
\frac
{1}{\ep}\sum_{n=0}^{\infty}\frac{B_{2n}}{(2n)!}\ep^{2n}
\int_0^x f^{(n)}(t)dt ,\ \ \ep \rightarrow 0 
\ee 
which is, of
course, the Euler-Maclaurin summation formula. Proposition \ref{P1}
below justifies our derivation of (\ref{eumclr}).

{\bf b.} Analogously, we can obtain easily the asymptotic
behavior of special functions from their generating 
recurrence relation.
For instance, the recurrence relation
\be{Besselrec}
 y_{k+1}+y_{k-1}=2(1+k\ep)y_k
\ee
has, for a fundamental set of solutions the Bessel functions 
$J_{k+\ep^{-1}}(\ep^{-1})$ and $Y_{k+\ep^{-1}}(\ep^{-1})$. 

Let $x=k\ep$, $\nu=k+\epi,\ a=\epi\nu^{-1}$. We fix $\rho<
\frac{2}{3}$
and for for $x>\ep^{\rho}$
 we take $y$ of the form (\ref{formalseries})
 determining the successive
terms by substituting the formal series in (\ref{Besselrec}).
To leading order
\be{'aprbess}
e^{\Phi'_0(x)}+e^{-{\Phi'_0(x)}}=2(1+x)
\ee
and we can choose two independent solutions
\be{expbes}
\epi \Phi_{0;\pm}=\pm \nu (\alpha -\tanh(\alpha ))\ \ (\alpha\equiv (\cosh )^{-1}(1/a))
\ee
For the next order $\Phi_{1;\pm}$ we obtain:
\be{bes2}
\Phi_{1;\pm}'/\Phi_{1;\pm}=-\frac{1}{2}\Phi_{0;\pm}''\coth(\Phi')
\ee
which can, again be integrated explicitly to give
$\exp(\Phi_{1;\pm})=const \sinh(\alpha)^{-1/2}$. That is, the asymptotic
behavior is 
\be{asybes}
(\nu \tanh(\alpha))^{-1/2}(A\exp[ \nu(\alpha-\tanh(\alpha))+\ldots ]+B\exp[-\nu(\alpha-\tanh(\alpha))+\ldots ])
\ee
and we get
 the familiar expressions in the theory of Bessel
functions (see [2], [3]). All the successive
orders are obtained easily in the same way.
For $x<-\ep^{\rho}$ we obtain from the same equations (\ref{'aprbess}), (\ref{bes2})
similar expressions but 
with trigonometric instead of hyperbolic functions. For small
values of $|x|$, $|x|~\  ~^<~\hskip -15pt ~_{\sim}~\ ~\ep^{2/3}$ 
the above asymptotic series becomes singular.
The appropriate  WKB-like series in this region  is in powers
of $\ep^{1/3}$ and the coefficients will be smooth
functions of $\ep^{1/3}k$, but otherwise the calculation
can be done explicitly in the same way; there is a region
of overlap where both asymptotic series are valid, namely
$\ep^{2/3}<|x|< \ep^{1/2}$ and where the series can be
matched.

 Proposition \ref{P2} can be used to make the above
approach rigorous.

{\bf c.} Finally, let $q$ be a smooth function and consider the Cauchy problem
\be{difex}
y''=q(x)y, \ \ y(0)=0,\ \ y'(0)=1
\ee
Assume for simplicity that 
$q:[0,1]\rightarrow \rdb^+$ and consider the associated
Euler scheme
\be{eulersch}
y_{k+\ep;\ep}+y_{k-\ep;\ep}=(2+\ep^2q(k\ep))y_{k;\ep},\ \ y_{0;\ep}=0,\ \ y_{1;\ep}=\ep
\ee
To characterize $y_{k;\ep} $ for small $\ep$ we substitute
\be{formalser}
y_{k;\ep}\ \sim\ \sum_{m=0}^{\infty}\ep^m \Xi_m(k\ep)
\ee
in (\ref{eulersch}). Note that if $\Xi_0\neq 0$ we obtain
this type of series from (\ref{formalseries}), when $\Phi_0=0$,
by expanding the exponential. 

The substitution leads to the equations:
\be{eqdifG}
\Xi''_m=q(x)\Xi_m-2\sum_{1\leq s;2s \leq m}\frac{\Xi_{m-2s}^{(2s)}}{(2s)!}
\ee
with the initial conditions
\be{inicond}
\Xi_m(0)=0;\ \ \ \Xi'_m(0)=-\sum_{k=1}^m \frac{\Xi_{m-k}^{(k+1)}(0)}{(k+1)!}
\ee
In particular, it follows that the scheme converges to the 
solution of the given Cauchy
problem (as it should) and because, as it is easily seen,  $\Xi_1\equiv 0$
the error $|y_{k+\ep}-\Xi_0(k\ep)|$ is $O(\ep^2)$. Proposition \ref{P3} applies to this example.

\section{Main results}
In this section we consider  the problem (\ref{mainrec}) and  the assumptions
following it  and give conditions under which the exact solutions
of the recurrence have asymptotic series to all orders in $\ep$.

\bp{P1} Let $\lambda_1(x),\ldots \lambda_l(x) $ be the roots of the
characteristic polynomial
\be{charpol}
\sum_{j=0}^l a_j(x,0)\lambda^j=0 
\ee
and assume that they are simple throughout I:
\be{distroo}
\inf_{x\in
 I}\{\mid\lambda_m(x)-\lambda_n(x)\mid\}>0\ \ {\mbox{if}}\  m\neq n.
\ee
(consequently, we will choose $\lambda_m$ 
to be $C^\infty$ functions).
 Suppose also that the interval $I$ is a finite union of (closed)
intervals such that in each one of them the ordering of the moduli
of the roots does not change, i.e.

$$
I=\bigcup_{p=1}^{P}I_p \ \ {\mbox{and}\  }\forall p\leq P \ \exists
(i_1\ldots  i_l) {\mbox{ such that}\ }
$$
\be{setS}
 |\lambda _{i_1}(x)|\leq |\lambda_{i_2}(x)|\leq \ldots \leq
|\lambda_{i_l}(x)| {\mbox{ for all}\ } x\in I_p.
\ee
where $(i_1\ldots  i_l)$ is a permutation of $(1\ldots  l)$.
Then there exists $\e'\leq \e_0$ and  l  functions
$F_1(x,\e)\ldots  F_l(x,\e)$ on $ I$ $\times [0,\e']$,
$C^\infty$ in each variable, such that for each
$\ep$,\break $\{\exp(\epi F_m(k\epsilon ,\epsilon ))\}_{m=1\ldots l}$ form a fundamental
set of solutions of the recurrence relation (\ref{mainrec}) in the sense that
for any solution $y_{k,\ep}$ of (\ref{mainrec}) there exist 
constants $C_1^{(p)},\ldots C_l^{(p)}$
such that in each $I_p$, 

\be{7} 
y_{k,\ep}=\sum_{m=1}^lC_{m}^{(p)}\exp[\epi F_m(k\epsilon ,\epsilon )]
\ee
\end{proposition}

\newtheorem{remark}{Remark}
\begin{remark}
In particular, this means that for small $\ep$
\be{8}
F_m(x,\ep)\sim \sum_{s=0}^{\infty}\Phi_{m,s}(x)\ep^s
\ee
where the $\Phi_{m,s}$ are smooth; they can be computed from
(\ref{mainrec}) by usual perturbation
expansions  in $\ep$. For example the first term is gotten from
the eikonal equation,

\be{9}
\exp(\Phi_{m,0}'(x))=\lambda_m(x)
\ee
giving the connection between the functions $F_m$ and the roots
$\lambda_m$ of
the characteristic polynomial. The second term is obtained from

\be{10}
\Phi_{n,1}'(x)=-\frac{\sum_{j=0}^l[j^2\Phi_{n,0}(x)''/2A_{j,0}(x)+\exp(\Phi_{j,0}'(x))A_{j,1}(x)]}{\sum_{j=0}^ljA_{j,0}(x)\exp(\Phi_{j,0}'(x))}
\ee
and so on.$\\$
\end{remark}
Note also that if not all the functions  $F_m$  have the same real
part, then there exist solutions of the recurrence which
are exponentially small relative to  the dominant ones and which
will therefore be unstable in the sense that a small ``generic''
perturbation of the initial condition will change completely
the behavior of these solutions; a perturbation series not involving
terms beyond all orders will only see the component of the 
solution along the dominant directions. However the relative size
of the solutions  could change with $x$ and then all the coefficients in 
the expansion (\ref{7}) could  be important for matching different regions.

\vskip 0.25cm
We now address the question of the asymptotic behavior of the
solutions of the recurrence when two characteristic roots cross. 
The previous Proposition is generalized below to the case when 
two of the roots
of the characteristic polynomial become equal at some point in $I$
provided  the
roots do not coalesce too quickly 
(condition (\ref{11'})). 
In this case a fundamental set of solutions has a more 
complicated structure.  Not too close
to the crossing point, a solution is still a linear combination of the 
form (\ref{7}) but the coefficients $C_m^{(p)}$ can change at the crossing (Stokes phenomenon). Very close of  the crossing, the coalescing roots bring in the
asymptotic expression of the solution.
series in noninteger powers of $\ep$.
This is in some sense the discrete counterpart of the
turning point-behavior of the solution of a differential equation
depending on a small parameter. 

Consider a subinterval $I_0$ of $I$ such that two roots
of the characteristic polynomial cross once in $I_0$ (say, at $x=0$)
and except for this, the ordering of the
moduli of the characteristic roots is constant
in $I_0$. The crossing is assumed to be generic:
\be{11'}
\left|\l_m(x)-\l_n(x)\right|>\ct\sqrt{|x|}\ ({\mbox{ for }}m\ne n)
\ee

To avoid excessive branching of the discussion and
formuli, we also assume
that the coalescing roots are {\it complex--conjugate}
for negative $x$ and {\it real--valued} for $x$ positive.
The general case is treated in a very similar way.

Fix two constants $\frac{1}{3}<\beta<\alpha<\frac{1}{2}$.
Then,

\bp{P2} 
i) For $|k|>\e^{-\beta}$
the general solution of the recurrence can be written
in the form
\be{regcase}
y_k=\sum_{j=1}^{l}C_j\exp(\ei F_j(k\e,\e))
\ee 
with $F_j$ as in Lemma \ref{LB}. The constants 
$C_j$ depend, in general on the sign of $k$.

ii)  For $|k|<\e^{-\alpha}$ a fundamental set of
solutions can be chosen such that:\break $l-2$ solutions are 
of the form $\exp(\ei F_j(k\e,\e))$ and two
special solutions of the form $$\exp(F_{\pm}(k\e^{1/3},\e^{1/3}))$$
where the functions $F_j$ and $F_{\pm}$  are smooth and 
$\exp(F_{\pm}(x,0))=$Ai$(\Theta\,x)\pm $Bi$(\Theta\,x)$
where Ai, Bi are the Airy functions (for the value of $\Theta$ see equation (\ref{Aisc}) .

Moreover there is a particular solution of the recurrence which 
has the behavior

$$y_k\sim Ai(\Theta\,k\e^{1/3})(1+\e^{1/3}A_1(\Theta\,k\e^{1/3})+..)$$
for large $k<\e^{-\alpha}$

The representations in (i) and (ii) are simultaneously valid in the
region $\e^{-\beta}<|k|<\e^{-\alpha}$ where the asymptotic
series can be matched.
\end{proposition}

{\bf Note}. A similar result can be proven if condition (\ref{11'}) is replaced
by
\be{16}
\left|\l_m(x)-\l_n(x)\right|>\ct\, x^{const'}
\ee
for some $\ct>0$.
Another case which is interesting for schemes that converge
to differential equations for small $\e$ is covered by the 
following proposition, in which we assume that the 
roots of the {\it complete} characteristic polynomial: 
\be{compchar}
\sum_{j=0}^l a_j(x,\ep)\lambda^j ,\ \ \ \ m=1\ldots l
\ee
are nondegenerate in a {\it higher} order in $\e$.

\bp{P3} Assume that the roots $\lambda_1(x,\ep)\ldots  \lambda_l(x,\ep)$
of (\ref{compchar}) satisfy the estimates
\be{degcond}
|\lambda_m(x,\ep)-1|=\ep^{q}(Q_{m}(x)+o(\ep)),\ \ m\leq l
\ee
for some  $q\in \ndd$ where the 
smooth functions $Q_{i}$ verify 
$$\inf_{x\in I}Q_{m}(x)>0,\ \inf_{x\in I}|Q_{m}(x)-Q_{n}(x)|>0
\ \ (m\neq n)\  m,n\leq l$$
Then, the conclusion of Proposition \ref{P1} holds and, furthermore
, for any formal series solution of the recurrence relation (\ref{mainrec})
there exists a true solution which is asymptotic to it.
\end{proposition}

In this particular case, since, as we shall see, $F_m(x,0)=0$, it
is more natural to represent the formal solutions as power series:

%%%%%%%%%%%%%%%%
\be{forso}
\Gamma=\sum_{i=0}^{\infty} \Psi_i(x)\ep^i
\ee
where the $\Psi_i$ are smooth and subject to the condition
%%%%%%%%%%%%%%
\be{locso}
\sum_{j=0}^l a_j(x,\ep)\Gamma_s(x+j\ep;\ep)=o(\e^s)\ \  \forall s,\ \ \ \ x\in I.
\ee
where $\Gamma_s(x;\epsilon)=\sum_{i=0}^s\Psi_i(x)\epsilon^i$.
The series (\ref{forso}) and those appearing
in the previous WKB expansions are usually divergent and one
could imagine that by iterating the recurrence 
the small error appearing in the local condition (\ref{locso})
could quickly reach  O(1).
Under the given restrictions, however,  Proposition \ref{P3} guarantees 
that $\Gamma_s$ is indeed an $o(\ep^s)$ approximation to a 
genuine solution.

It can be now checked without difficulty that Propositions \ref{P1}, \ref{P2} and \ref{P3}
apply to the examples a), b) and c), respectively.

The layout of the paper is as follows: in Section~3 we prove
our main results and in Section~4 we discuss some further 
applications of these results.

\section{Proof of the results} To prove Proposition \ref{P1} we show (Lemma \ref{LA}) 
the existence of $l$ formal series solutions of the form (\ref{formalseries}) to
the recurrence. 

We then show (Lemma \ref{LB}) that the the proposition is true
if $P=1$ (cf.(\ref{setS})). The proof is by induction on the
order $l$ of the recurrence. First, we choose the particular
formal solution corresponding (cf. (\ref{18})) to the root with maximum modulus 
---which gives the ``stable'' direction--- and show that there is
a true solution with this asymptotic behavior. Next we use this particular
solution to decrease the order of the recurrence by 1.

{\bf 3.1} The functions
$\Phi_{m,s}$ of the following lemma will turn out to be the 
functions giving the asymptotic expansions (\ref{8}) and can
be obtained by requiring that (\ref{mainrec}) is
satisfied in all orders in $\e$ by the formal solution
$y_m=\exp(\ei\sum_s\Phi_{m,s}(x)\e^s)$.

\bl{LA} For each $m=1\ldots l$ there exists a sequence $\{\Phi_{m,s}\}_{s\in 
{\Nsm }}$ of functions in $ C^\infty ({ I})$ such that 
\be{17}
\exp(-\epi\Phi_{m,0}(k\epsilon ))\sum_{j=0}^l a_j(k\epsilon ,\epsilon
)\exp(
\epi\
\sum_{t=0}^s \ep^{t}\Phi_{m,t}((k+j)\epsilon )))={\mbox{O}}(\ep^{s+1})
\ee
and
\be{18}
\exp(\Phi_{m,0}'(x))=\lambda_m(x)
\ee
for $\ep \leq \e_0$, 
$s \in \ndd$ and $k\epsilon  \in { I}$
\end{lemma} 

The proof of Lemma 1 is by induction on the expansion order $s.$
In view of (\ref{asona}) and (\ref{asnonz}) we can define for each characteristic
root $\lambda_m(x)$ (i.e. root of eq. (\ref{charpol})),  a function  $\Phi_{m,0} \in C^\infty(I)$
such that (\ref{18}) holds.

It is then straightforward to show that
$\exp(\epi\Phi_{m,0})$ verifies (\ref{mainrec}) to ${\mbox{O}}(\epsilon )\exp(\epi\Phi_{m,0}(x))$ so that (\ref{17})
holds for $s=0$.
Assuming now that $\Phi_{m,0},\Phi_{m,1}\ldots  \Phi_{m,s_0}$ are already
defined so that for all $s\leq s_0$ (\ref{17}) is verified, one can easily
check that for any $\Psi\in C^\infty ({ I})$,
%%%%%%%%%%%%%%%%%%%%%%%%%%%%%%%%%%%%%%%%%%%%%%%%%
$$\exp[-\epi\Phi_{m,0}(k\epsilon )]$$
%%%%%%%%%%%%%%%%%%%%%%%%%%%%%%%%%%%%%%%%%%%%%%%%%
$$\sum_{j=0}^l \exp\{\epi\sum_{t=0}^{s_0}\Phi_{m,t}((k+j)\epsilon)\ep^{t}+\ep^{s_0+1} \Psi((k+j)\epsilon )\}a_j(k\epsilon ,\epsilon )=$$
%%%%%%%%%%%%%%%%%%%%%%%%%%%%%%%%%%%%%%%%%%%%%%%%%%%%%%%%%%%%%%%%%%%%%%%%%%%
\be{19}\ep^{s_0+1}[\Psi '(k\epsilon )\sum_{j=0}^l j\exp\{j\Phi '_{m,o}(k\epsilon )\}
a_j(k\epsilon ,0)
+H_{s_0}(k\epsilon ,\epsilon )]+{\mbox{O}}(\ep^{s_0+2})
\ee
where $H_{s_0}$ is a smooth function.

Since by (\ref{asnonz}) and (\ref{distroo}) 
\be{20}
\inf_{x \in { I}} |\sum_{j=0}^l j
a_j(x,0)\lambda_m^j(x)|>0
\ee
one can define a smooth function $\Psi (x)\equiv\Phi_{s_0+1}(x)$ 
such that the term
in square brackets in the RHS of (\ref{19}) vanishes.$\Box_{L1}$

We note at this point that the series $\sum_{s=0}^\infty
\Phi_{m,s}(x)\ep ^{s}$ is, usually, not convergent and  there does
not yet follow the existence of a solution asymptotic to it.

Now we address the question of existence of true
solutions of the recurrence having the precribed
asymptotic behavior.
In what follows, we shall understand by a formal  solution
an expression $\tilde Y=\exp(\ei\sum_{s=0}^{\infty}\Phi_{m,s}(x))\e^s$ 
satisfying the conclusions of Lemma \ref{LA}. Given 
$\Phi_{m,0}$, the $\Phi_{m,s}$ are uniquely determined up to
integration constants. 

Assume first $P=1$ (cf. (\ref{setS})).
Relabelling if necessary, we  assume that for $m\leq n$,  $|\lambda_m(x)|\leq |\lambda_n(x)|$ on 
$I$.

\bl{LB}Let  $\tilde Y$  be a formal  solution of
(\ref{mainrec}) ($m \in\{1\ldots l\}$ fixed). 
 There exists a sequence
$\{Y_{m,k;\ep}\}_{k,\ep}$ such that for any $\ep\leq\e_0$, $Y_{m,k;\ep}$ is a solution of 
(\ref{mainrec}) for $k\epsilon  \in I$ having $\Gamma_m$ as an asymptotic series in the
sense that there is a sequence of positive constants $\{C_s\}_s$ such that 
\be{21}
|Y_{m,k;\ep}\exp[-\epi\sum_{t=0}^s\Phi_{m,t}(k\epsilon )\ep^{t}]-1|<C_s\ 
\ep^{s+1} \ for \ all\ s\in\ndd
\ee
\end{lemma}
\begin{remark}
The previous lemma can be restated as follows: for each m=1\ldots l
 there is a function
$F_m(x,\e):I \times [0,\e_0] \rightarrow {\cdb}$, smooth in each
variable, such that $exp(\epi F_m(k\epsilon ,\epsilon ))$ is a solution of the
recurrence (\ref{mainrec}) for every $\epsilon$, $k\epsilon$  $\in I$ and, as $\e\rightarrow 0$,
\be{22}
F_m(x,\e)\sim \sum_{t=0}^\infty \Phi_{m,t}(x)\e^t
\ee
\end{remark}
This remark follows easily from a classical result ( see e.g. 
[4], page 33 and [5]) 
stating that for any sequence of
numbers there is a smooth function having that sequence for its 
derivatives at the origin and from the easily proven fact that for 
each sequence
$\{(x_n,a_n)\}_{n\in{\Nsm}}\ ,\ x_n \searrow 0 $ and $a_n n^k \rightarrow 0$ $ \forall k$ one 
can construct a $C^\infty$ ``interpolation'' function f such that
 $f^{(k)}(0)=0$ $
\forall k$ and $ f(x_n)=a_n$ for all $n$.

Some of the estimates that we need for proving Lemma 2 
become to a certain extent easier by the remark
that a global shift in the estimates defining an  asymptotic
series is unimportant:
%#####################
\begin{remark}

Let F be a function such that for a fixed $s_0\geq 0$ and any s in $\ndd$
\be{23}
\limsup_{\e\rightarrow0}\e^{-s-1}|F(\e)-\sum_{t=0}^{s+s_0}C_t \e^t |=D_s
\ee

Then
\be{24} 
\limsup_{\e\rightarrow0}\e^{-s-1}|F(\e)-\sum_{t=0}^{s}C_t \e^t |\leq D_s +
|C_{s+1}|
\ee

\end{remark}
%#########################
Thus,  to prove (\ref{21}) we need only  show that
for some fixed ${s_0}$, (we drop the subscript $\e$
to ease the notation)
\be{21'}
|Y_{m,k}\exp(-\epi \sum_{t=0}^s\Phi_{m,t}(k\epsilon
)\ep^{t})-1|<C_s\ 
\ep^{s+1-s_0}
\ee

\begin{remark}
In view of condition (\ref{distroo}) it is easy to check that, for
small $\ep$, the solutions 
$\exp(\epi F_m(k\epsilon ,\epsilon ))$ are linearly independent, thus forming a
fundamental set of solutions on $I$.
\end{remark}

{{Proof of Lemma \ref{LB}}} The proof is by induction on the order $l$ of the recurrence.
\vskip 0.25cm
\newtheorem{step}{Step}
\begin{step}
For any l, in the same conditions as in Lemma \ref{LA}, the estimate
(\ref{21}) holds for m=1 {\mbox {(recall that $\lambda_1(x)$ is the largest in
absolute value)}}.
\end{step}
First we show that if  a solution is asymptotic to the formal solution
corresponding to $\lambda_1$
at the left end of $I$ it remains asymptotic to it throughout $I$.
We can assume without loss of generality
that the left end of $I$ is at  $x=0$; at this point we   choose 
appropriate initial conditions: 

let  $\eta_p,\ p=1\ldots l$ be any functions such that as $\e\rightarrow 0$
\be{25}
\eta_p(\e)\sim \exp(\e^{-1}\sum_{t=0}^\infty \Phi_{1,t}(p\e)\e^t)
\ee
(see Remark 2) and let 
\be{26}
\tilde{Y}_{\ep;s}(x)=\exp(\epi\sum_{t=0}^s \Phi_{1,t}(x)\ep^{t})
\ee

Let also $\ \ Y_{1,k}\ \ $ be the\ \  solution of \ (\ref{mainrec})\
\   satisfying
the\ \  initial conditions \\
$Y_{1,p}=\eta_p(\epsilon )\ ,\ p=1\ldots l$. It is natural to
rescale the recurrence relative to its approximate solution: let
\be{27}
C_{k ;s}=Y_{1,k}/\tilde{Y}_{\ep;s}(k\epsilon )
\ee
Then the recurrence relation for $C_{k ;s}$ can be written 
%%%%%
\be{28}
\sum_{j=0}^l \tilde a_{j,s}(k\epsilon ,\epsilon )C_{k+j ;s}=0
\ee
%%%%%
where $ \tilde{a}_{j,s}(k\epsilon ,\epsilon )={a}_j(k\epsilon ,\epsilon
)\tilde{Y}_{\ep ;s}((k+j)\epsilon )/
\tilde{Y}_{\ep ;s}(k\epsilon )$. It can be seen that (\ref{28}) is of the same type
as (\ref{mainrec}) and that, for $\ep$ small enough, it satisfies the corresponding
assumptions (\ref{asnonz}), (\ref{distroo}) and (\ref{setS}) on $I$.

Also,  
$\max\{|\tilde\lambda_m(x)|;x\in I,\ m=1\ldots l\}=1$.
Now, (\ref{21'}) means  that  for some 
fixed $s_0$,
\be{eqcond}
|C_{k ;s}-1|<const_s\ep ^{s-s_0+1}
\ee
From the definition of $\tilde {Y}_{\ep;s}$ it follows that
\be{approxid}
\sum_{j=0}^l \tilde a_{j,s}(k\epsilon ,\epsilon )=O(\ep^s)
\ee
Rewriting the recursion relation \ (\ref{28})\  in
matrix form\ \ \  ${\bf C}_{k+1}~=~\tilde A_{k}{\bf C}_k $, where
%%%%%%%%%%%%%%
\be{29}
\tilde A_{k}=\left(
\begin{array}{cccc}
 \frac{-\tilde a_{l-1}(k\epsilon ,\epsilon )}{\tilde a_l(k\epsilon ,\epsilon )} & 

\frac{-\tilde a_{l-2}(k\epsilon ,\epsilon )}{\tilde a_l(k\epsilon ,\epsilon )} &... &  \frac{-\tilde a_{0}(k\epsilon ,\epsilon )}{\tilde a_l(k\epsilon ,\epsilon )} \\

1 & 0 & ... & 0 \\ 
0 & 1 & ...& 0 \\ ... & ... & ... & ...\\ 0 & 0 & ..1 & 0 \end{array} 
\right)
\ee
and rewriting also (\ref{approxid}) as
\be{approxidm}
\tilde A_{k} {\bf 1}={\bf 1}+\ep^{s}{\bf E}_{k}
\ee
where ${\bf 1}_j=1$ and $\|{\bf E}_{k}\|<const$ uniformly in $k,\ep$ we have
\be{explicsl}
{\bf C}_k={\bf 1}+\ep^s\sum_{j=1}^k \tilde A_{k}
\tilde A_{k-1}\ldots \tilde A_{j+1}{\bf E}_{k}
\ee
Step 1 is completed by showing (\ref{eqcond}) (and thus (\ref{21'})),
which follows from the stability lemma below.
\bl{LC} Let $\tilde A_{k}$ be a family of matrices of the form (\ref{29})
where  $\tilde a_j:I\times [0,\e_0]\rightarrow{\cdb}$ (I is an interval)
as in Lemma \ref{LB}.

Suppose further that the roots $\tilde \lambda_m(k\ep,\ep)$ of the
polynomial
%%%%%%%%%%%%%%%%%
\be{complepol}
\sum_{j=0}^l \tilde a_j(k\ep,\ep)\tilde\lambda^j=0 ,\ \ \ \ m=1\ldots l
\ee 
%%%%%%%%%%%%%%%% 
satisfy
\be{30}
\sup_{k\ep\in I}\{|\tilde\lambda_m(k\ep,\ep)|\}\leq 1+const\ \ep,\ m=1..l
\ee
%%%%%%%%%
and that the condition corresponding to (\ref{distroo}) is fulfilled on I.
Then there is an $\ep$-independent constant C such that

\be{31}
\|\tilde A_{k}\tilde A_{k-1}\ldots \tilde A_{j+1}\|\leq 1+C |k-j|\ep
\ee

\end{lemma}

{ {Proof of Lemma 3.}}
The eigenvalues of the matrix
(\ref{29}) are the $\tilde\lambda_m(k\ep,\ep)$ and the corresponding
eigenvectors matrix is  $(\tilde{\Lambda}_{k})_{i,j}~=~
\tilde\lambda_j(k\ep,\ep)^{l-i}$. We then write the product on the 
LHS of (\ref{31}) as

\be{32}
\tilde \Lambda_{k}D_{k}\tilde \Lambda_{k}^{-1}\tilde
 \Lambda_{k-1}D_{k-1}\ldots \tilde
 \Lambda_{j+1}D_{j+1}\tilde \Lambda_{j+1}^{-1}
\ee

where
\be{33}
D_{p}= diag(\{\tilde\lambda_m(p\epsilon ,\epsilon )\}_{m=1\ldots l})
\ee

and the proof follows from (\ref{30}) and the estimate
\be{34.0}
\|\tilde \Lambda_{p}^{-1}\tilde \Lambda_{p-1}\|\leq 1+const\ \epsilon 
\ee
which can be checked, for instance,  using the following explicit formula,
whose elementary proof we omit:
\begin{remark}
Let X and Y be two nonsingular Vandermonde-type matrices $X_{i,j}=x_j^{l-i},
\ Y_{i,j}=y_j^{l-i},\ i,j=1\ldots l$. Then,
\be{Vandermonde}
(X^{-1}Y)_{i,j}=\prod_{n\neq i}\frac{y_j-x_n}{x_i-x_n}
\ee
\end{remark}
\vskip 0.25cm
\begin{step}\             
\end{step}
The conclusion of Lemma \ref{LB} for $l=1$ follows from step 1.

 Now we assume that the conclusion of the lemma holds
for all recurrences of order less than $l$ and
prove it for order $l$, by reduction to the
$l-1$ case. In view of the first step, we know
that to $|\l_1|$ there corresponds a true solution $Y_1$
for which the asymptotic behavior is the formal
solution $\tilde Y_1 $. We shall use this solution
to reduce the order of the recurrence by one. Let
\be{35} 
C_{k}=y_{k}/Y_{1,k}
\ee
The recurrence relation for $C_k$  is then of the form (\ref{28}) where now
\be{36}
 \tilde{a}_j(k\epsilon ,\epsilon )={a}_j(k\epsilon ,\epsilon ){Y}_{1,k+j}/
{Y}_{1,k}
\ee
and obviously, the asymptotic behavior to all orders is the same as if we had
made the rescaling with respect to a formal solution.
The point is that now, instead of (\ref{approxid}) we have
\be{37}
\sum_{j=0}^l\tilde a_j(k\epsilon ,\epsilon )=0
\ee
so that the  $y=1$ is an actual solution. To use this fact, let
$d_{k}=C_{k+1}-C_{k}$. We get,

\be{38}
\sum_{s=0}^{l-1}b_s(k\epsilon ,\epsilon )d_{k+s}=0
\ee
where $b_s(k\epsilon ,\epsilon )=\sum_{j=s+1}^{l}\tilde a_j(k\epsilon ,\epsilon )$.

The
characteristic equation for (\ref{38}) can be written as

$$\sum_{j=1}^l \tilde a_j(x,0)\sum_{s=0}^{j-1}\tilde\lambda^s=0$$
or, for $\lambda\neq 1$, as it easily follows from (\ref{37}),
\be{39}
\sum_{j=0}^l\tilde a_j(x,0)\tilde\lambda^j =0
\ee
We first  check that the new recurrence satisfies the
hypothesis of the Lemma. But this
is easy since the new coefficients 
are finite combinations of $a_j$ and in particular
$b_0=\sum_{j=1}^{l}\tilde a_j=-a_0$ and $b_{l-1}=a_l$,
and in view of (\ref{39}) the same arguments as in
step one apply to see that the characteristic roots
have the required properties. Then we want to check that
we have the required number of appropriate formal solutions.
This is also straightforward because we can derive
them from the formal solutions of the original equation.
Indeed, 

\be{defd}
d_k=y_{k+1}/Y_{1,k+1}-y_{k}/Y_{1,k}
\ee
If we substitute for  $y$ a {\it formal} solution $\tilde Y_m$
we obtain the formal expression : 

$$\tilde d_k=\exp\{\epi\sum_{t=0}^{\infty}[\Phi_{m,t}((k+1)\epsilon
)-\Phi_{1,t}((k+1)\epsilon )]\ep ^{t}\}$$

\be{41}
-\exp\{\epi\sum_{t=0}^{\infty}[\Phi_{m,t}(k\epsilon
))-\Phi_{1,t}(k\epsilon )]\ep ^{t}\}=
\ee
$$\exp\left\{\e^{-1}(\Phi_m(x)-\Phi_1(x))
+{\mbox{ series}}\right\}[\l_m(x)/\l_1(x)](1+\epsilon\times{\rm series}]$$

which, since $\l_m(x)/\l_1(x)$ does not vanish
can be written 
 as an exponential of a formal series, the form required 
by our arguments:
\be{42}
\exp(\epi\sum_{t=0}^{\infty}\Delta_{m,t}(k\epsilon ))
\ee
and then the expressions (\ref{42}) for $m=2\ldots l$ are formal
solutions for the recurrence of order $l-1$ (\ref{38})
---because, by construction, (\ref{41}) are; it follows, by the 
induction hypothesis that there exist  true solutions of (\ref{38})
of the form
\be{43}
d_{m,k}=\exp(\epi D_m(k\epsilon ,\epsilon ))
\ee
where $ D_m(\cdot,\cdot)$ are smooth functions
having the asymptotic behavior given by (\ref{42}).
\vskip 0.50cm
\begin{step}\ 
\end{step}
 To complete the
proof of Lemma 2 it remains only to check that 

\be{44}
\exp\left(\epi F_1(k\epsilon ,\epsilon )\right)
               \sum_{p=0}^k \exp(\epi D_m(p\epsilon ,\epsilon ))
\ee

\z has the asymptotic behavior needed for the original recurrence, i.e.

\be{45}
\exp(\epi\sum_{t=0}^\infty(\Phi_{m,t}(k\epsilon )\ep ^{t}))
\ee
We let $k_1$ ($k_2$)  be the left (right, respectively) end of the
interval. Both $k_1$ and $k_2$ might depend on $\e$.
 By the definition of the
$d_{k;m}$ we have 

$$C_{m;k}=C_{m;k_1}+\sum_{i=k_1}^k d_{m;i}$$

With the choice 

$$C_{m;k_1}=-\sum_{i=k_1}^{k_2}d_{m;i}$$
we get the particular solution

$$C_{m;k}=\sum_{i=k}^{k_2}d_{m;i} $$
from which, referring to the definition of the 
$C_{m;k}$, we get a solution of the $l-1$ recurrence
in the form (the choice of the sign will
become clear later)

$$Y_{m;k}=-Y_{1;k}\sum_{i=k}^{k_2}d_{m;i} $$
whose the asymptotic behavior is given by the
formal solution $\tilde Y_{m}$. Indeed,
for any  fixed large $s$ we have,
$$Y_{m;k}=-Y_{1;k}\left[\sum_{i=k}^{k_2}
                  \exp\left(\ei\sum_{t=0}^{s}\Delta_t(i\e)\e^t\right)
                  \left(1+O(\e^{s-1})\right)\right]=$$
$$-Y_{1;k}\left[\sum_{i=k}^{k_2}\left(
                  \exp\left(\ei\sum_{t=0}^{s}\Phi_{m;t}(i\e+\e)\e^t
                   - \Phi_{1;t}(i\e+\e)\e^t\right)-
               \right.\right. $$
$$\left.\left.\exp\left(\ei\sum_{t=0}^{s}\Phi_{m;t}(i\e)\e^t
                   - \Phi_{1;t}(i\e)\e^t\right)
                  \right)\left(1+O(\e^{s-1})\right)\right]=$$

$$Y_{1;k}\left[\exp\left(\ei\sum_{t=0}^{s}\Phi_{m;t}(k\e)\e^t
                   - \Phi_{1;t}(k\e)\e^t\right)-
           \exp\left(\ei\sum_{t=0}^{s}\Phi_{m;t}({k_2})
                   - \Phi_{1;t}(k_2)\e^t\right)\right]+$$
\be{finest}
({k_2}-k)\max_k\left|
                 \exp\left(\ei\sum_{t=0}^{s}\Delta_t(i\e)\e^t\right)\right|
                  O(\e^{s-1})
\ee

Because, by assumption, $\lambda_1$ has
the largest modulus,
$\Re\left(\Phi_{m;0}(k\e)-\Phi_{1;0}(k\e)\right)$ 
is nonincreasing in $k$ in the given region. Therefore
 (\ref{finest}) equals,

$$Y_{m;k}\left(1+\ct\,\max_k\left\{\left|\exp\left(\Phi_{m;1}(k\e)
                   - \Phi_{1;1}(k\e)\right)\right|\right\}o(\e^{s-2})\right)=$$
\be{endest} 
Y_{m;k}\left(1+o(\e^{s-2})\right) 
\ee

  \hskip 10cm                                       $\Box_{L2}$

The proof of Proposition \ref{P2} follows essentially the same steps 
but is, as expected, more involved in the regions
of near--breakdown of the asymptotic series. The details
are given in the next section.

The proof of Proposition \ref{P3} is very easy, using an estimate of the
form (\ref{31}) for the matrices corresponding to the original recurrence,
estimate which is straightforward to obtain from Remark 5 and the
hypothesis of the proposition.

\begin{section}{Proof of Proposition 2.2}

We assume at first that the crossing occurs between the two largest
characteristic roots and explain at the end of the
proof how the general case is reduced to this one.

The layout is as follows.  We first study the small region around 
the crossing point (the interior region)
 where  $\exp(\ei\sum\Phi_{1,2;t}(x)\e^t)$ fail
to be  formal solutions (and the series occuring at the exponent
cease to be  asymptotic series). The new formal solutions
are to leading order combinations of Airy functions. Their formal
properties (domain of asymptoticity, growth in $x$) 
are examined. Next we show that there exist true
solutions of the recurrence that are asymptotic to them.
It is also shown that there exists a particular true solution  
which is asymptotic, to leading order, to the function
Ai($x\e^{-2/3}$) and which is important
for the matching problem (it gives the exponentially
decaying formal solution). 

We then show that the formal solutions coming
from the exterior region continue to represent
correctly the solutions of the recurrence 
far enough into the interior region (down to $|x|\sim\e^{2/3}$)
to allow for matching with the interior ones, 
which are valid up to $|x|\sim\e^{1/2}$. 

{\bf A. The interior region of the crossing interval}

This is the region $|x|\ll\e^{\frac{1}{2}}$; for definiteness
we fix $\alpha\in(\frac{1}{2},\frac{2}{3})$  and take it to be

$$D_{\a}= \{k:|k|<\epsilon^{-\alpha}\}$$
or, in terms of $\xi:=k\e^{1/3}:=k\d$ which is, 
for reasons that will become clear later,
the natural variable in this region,

\be{unbdo}
D_{\a}= \{\xi:|\xi|<\d^{1-3\alpha}\}
\ee

The basic steps of the proof of existence of solutions with given
asymptotics
are the same as for Lemma \ref{LB}. We will first 
obtain a solution corresponding to (one of the two) 
largest eigenvalues and with it reduce the problem
to a lower order, nondegenerate recurrence.

Because in the variable $\xi$ $D_{\a}$ is  unbounded,
a slight extension of Lemma \ref{LB} is needed for the interior
region. We now allow the interval $I$ in Lemma \ref{LB} to be 
of the form (\ref{unbdo})  but  strengthen the other
hypothesis. In order to make the correspondence with
Lemma \ref{LB}, note that $\xi$ plays the role of $x$
and $\d$ is the counterpart of $\e$. We require
the same conditions as in Lemma \ref{LB} and in addition,

 {\bf a)} $a_j(\xi,\d)$ are assumed to 
have asymptotic series $\sum_t\,a_{j,t}(\xi)\d^t$ valid
throughout  the region $D_{\a}$ which are smooth 
in the sense that the all their formal derivatives with respect to
$\xi$, 
$\sum_t\,a^{(m)}_{j,t}(\xi)\d^t$,  exist and $|a^{(m)}_{j,t}(\xi)|<\ct\,\d^{const'\,t}$

 {\bf b)} The roots of the characteristic polynomial 
\be{charp}
\sum_{j=0}^l a_j(\xi,0)\l(\xi)^j=0
\ee
are nondegenerate:$ $
\be{separ} {\inf}_{D_{\a}}\left|\lambda_m(\xi)
-\lambda_n(\xi)\right|>{\mbox{const}}>0\ (m\ne n)
\ee
and the polynomial itself is nondegenerate in the sense
\be{nondpoly}{
\inf}_{D_{\a}}\left\{|a_0(k\d)|,|a_l(k\d)|,\frac{1}{|a_j(k\d)|}\right\}>{\mbox {const}}>0
\ee

{\bf c)}
$$ \left|{\l_m((k+1)\d,\d)-\l_m(k\d,\d)}
                                \right|\,
                                      <\,\frac{{\mbox{const}}}{|k|+1}$$
in $D_{\alpha} $ where $\l_m((k+1)\d,\d)$ are the roots of the complete polynomial
$\sum_{j=0}^l a_j(\xi,\d)\l^j$

\bl{A1}
Under these assumptions
for any formal solution of (\ref{mainrec}) of the form

$$S:=\exp(\d^{-1}\Phi_0(\xi)+\sum_{m=0}^{\infty}\Phi_m(\xi)\,\d^m)$$
where the exponent is assumed to be a smooth asymptotic series
in the sense defined in a), there exists a true solution
of the recurrence which is asymptotic to it  in $D_a$.

\end{lemma}

{\bf Proof}. The proof follows closely the proof of 
Lemma \ref{LB}. We only emphasize the differences:
For the recurrence (\ref{28}) we have also to verify
condition c) :

$$\left|\tl_m((k+1)\d,\d)-\tl_m(k\d,\d)\right|=\left|\frac{\l_m((k+1)\d,\d)}{\l_1((k+1)\d,\d
)}-\frac{\l_m(k\d,\d)}{\l_1(k\d,\d)}\right|\,<$$
$${\mbox{const}\,}\left(|k|+1\right)^{-1}$$
The equivalent of Lemma \ref{LC} states now that there is
an $\d$-independent constant C such that

\be{matunb}
\| A_{k} A_{k-1}\ldots  A_{j+1}\|\leq (1+C
|k-j|^{{\mbox{const}}})\ \ \ 
\ee

Indeed, by  the Remark 5 a diagonal term of the matrix  
$T:=\Lambda_{p}^{-1}
\Lambda_{p-1}$ is of the form 
$$T_{m,m}=\prod_{n\ne m}\frac{\l_m((p-1)\d,\d)-\l_n(p\d,\d)}
{\l_m(p\d,\d)-\l_n(p\d,\d)}=
$$

\be{estdia}=\prod_{n\ne m}\left(1+\frac{\l_m((p-1)\d,\d)-{\l_m(p\d,\d)}}
{\l_m(p\d,\d)-\l_n(p\d,\d)}\right)=1+O((|k|+1)^{-1})
\ee
by c). Similarly, the moduli of the nondiagonal
terms are seen to be less than ${\mbox{const}}/(|k|+1)$. Therefore $T={\bf I}+R$ where the $\|R\|$ is 
 $O((|k|+1)^{-1})$, 
hence the inequality (\ref{matunb})  follows.

The last part of the proof of Lemma \ref{LC}
applies here without any significant change.\ \ \ $\Box_{L4.1}$
\vskip 0.25cm

\z The reduction to the nondegenerate case.
\vskip 0.25cm
\z We consider the initial recurrence in the 
neighborhood of a crossing point, say x=0  where
$\l_1(0)=\l_2(0)=1$ (the value at zero can be
chosen through a trivial global rescaling
of the recurrence).
It is convenient to consider rescaled variables 
$\d=\e^{1/3}$ and $\xi=k\d$. In these variables,
the coefficients $a_j(x,\e)=
a_j(\xi\d^2,\d^3)$   have smooth
asymptotic series in $\d$ in $D_{\alpha}$ which
are in fact obtained through series expansion in $x$
from (\ref{asona}):

\be{poly}
a_j(\xi\d^2,\d^3)\sim\sum_{s\ge 0}P_{j,s}(\xi)\d^{s}
\ee
where 
$$P_{j,0}=a_j(0,0);\ \ \ P_{j,1}(\xi)=0;\ \
P_{j,2}(\xi)=(D_x a_j)(0,0)\xi;
$$
and in general $P_{j,s}(\xi)$ are polynomials in $\xi$ of degree
at most $s/2$ for $s$ even and $(s-3)/2$ if $s$ is odd.
 To avoid complicating the  
notation we write $a_j(\xi,\d)\equiv a_j(\xi\d^2,\d^3)$.
We have first to find formal solutions for this new
recurrence.
\bl{Aux}
There exist $l$ linearly independent
formal solutions in $D_{\alpha}$, of the form 

\be{expform}
{\exp\left({\d^{-1}\sum _{t=0}^{\infty
}\Psi_{m,t}(\xi)\delta^{t}}\right)}
\ee

where $\Psi_{m,t}(\xi)$ are smooth in $\xi$ and satisfy 
the estimates

$$|\Psi_{m,t}(\xi)|<\ct_{m,t}+\ct'_{m,t}|\xi|^{\frac{t}{2}+1}$$

\end{lemma}

This means in particular that the domain of formal validity of the power series is 
then $\xi^{1/2}\d\ll1$ i.e., $x\ll1$. The domain
in which it is actually asymptotic to the solution
is however much smaller ($x\ll\sqrt\e$) as we shall see.

Proof of Lemma 4.2.

 The formal solutions corresponding to the nondegenerate
roots give rise automatically  to acceptable formal solutions
in the new variables $\xi,\d$. Indeed,
\be{ai1}
{\exp\left({{{\epsilon}^{-1}{\sum _{t=0}^{\infty }\Phi_{m,t}(x)\epsilon^{
t}}}}\right)}=
\ee

$${\exp\left({\d^{-3}\sum _{t=0}^{\infty }\Phi_{m,t}(\xi\,\delta^{2})
\delta^{3\,t}}\right)}=$$
$${\exp\left({\d^{-3}\Phi_{m,t}(0)+\d^{-1}\sum _{s,t\ge 0}^{\infty }\Phi_{m,t}^{(s)}(0)
\xi^s\delta^{2\,s+3\,t-2}}\right)}$$
The term in $\d^{-3}$ is merely a multiplicative constant so it 
can be dropped and we are left with
a formal solution of the form

\be{expfrm}
{\exp\left({\d^{-1}\sum _{t=0}^{\infty
}\chi_{m,t}(\xi)\delta^{t}}\right)}
\ee
where the $\chi_{m,t}(\xi)$ are in fact polynomials in $\xi$ 
of degree $\le t/2+1$.

For $m=1,2$ it is more convenient to write first the possible
formal series solutions for the equation, in the form:

\be{linform}
\sum _{t=0}^{\infty }\chi_{t}(\xi)\delta^{t}
\ee
and then show that we can write them 
in the form (\ref{expfrm}).

 Substituting (\ref{linform})  in the recurrence we get
$$\sum _{j=0}^{l}a_{j,\xi,\delta}\sum _{t=0}^{\infty }\chi_{t}
(\xi+j\delta)\delta^{t}$$
The term of order $s$ in (\ref{linform}) is  gotten by
differentiating the auxilliary
 equation 

$$\sum_{j=0}^{l}a_j(\xi,\delta)\chi(\xi+j\d,\d)=0$$
$s$ times with respect to $\d$.
$s$ times with respect to $\d$.
We get (see (\ref{poly}) )
$$\left.\sum_{j=0}^{l}\sum_{t=0}^{s}P_{s-t}(\xi)\left(D_{\xi}+
D_{\d}\right)^t\chi\right|_{delta=0}=0$$
which after expansion, change of order of summation and use of
(\ref{poly}) gives,
$$\sum_{j=0}^l\left[\frac{1}{{{s\choose s-2}}}\sum_{\sigma=0}^{s-3}
                             \sum_{t=\sigma}^{s}\left(
                      {s\choose t}{t\choose \sigma}P_{s-t}(\xi)
                      j^{t-\sigma}D_{\xi}^{t-\sigma}\chi_\sigma(\xi)\right)+\right.$$
\be{biginduct}
\left.
                     D_x a_j(0,0)\xi\,\chi_{s-2}+
                      a_j(0,0)j^2\chi_{s-2}\right]=0
\ee

It follows  that $\chi_0$ is obtained as a solution 
of the homogeneous Airy equation
\be{susu}
\chi''(\xi)=\Theta^3
\xi\,\chi(\xi) 
\ee
where 
\be{Aisc}
\Theta^3=\frac{\sum_{j=0}^{l}D_x
a_{j}(0,0)}{\sum_{j=0}^{l}a_{j}(0,0)\,j^2}
\ee

\z and that, given $\chi_0...\chi_{s-3}$,  we get $\chi_{s-2}$
as a solution of an inhomogeneous Airy equation of the form

$$\chi''(\xi)=\Theta^3
\xi\,\chi(\xi)+R(\xi)$$
where $R(\xi)$ is a linear combination of higher derivatives
of $\chi_0..\chi_{s-3}$. To avoid cumbersome notations,
we shall assume in the following that $\Theta$ is one.
We can  check that the assumption of genericity
($|\l_1(x)-\l_2(x)|>\ct\,\sqrt{x}$) implies  
$\sum_{j=0}^{l}D_x a_{j}(0,0)\ne 0$. We shall assume
for definiteness that it is negative. 

It follows by an obvious induction that the
$\chi_s$ are smooth. Now we show that they satisfy
the inequalities stated in the Lemma 4.2. 

${\underline{\rm Remark}}$ Consider the inhomogeneous Airy equation
$f''(\xi)=\xi\,f(\xi)+R(\xi)$ and assume 
 $R(\xi)\sim\xi^{\rho}\exp(2A/3\xi^{3/2})$ with
$A=\pm 1$ for $x\rightarrow\infty$ and $A=\pm i$ at
$-\infty$ . Then $f(\xi)\sim$ $\xi^{\rho-1}\exp(2A/3\xi^{3/2})$.
This estimate follows immediately from the explicit form of the solution:

$${f(\xi)=\mbox{
Ai}(\xi)\int^{\xi} R(t)Bi(t)\,dt-Bi(\xi)\int^{\xi} R(t)Ai(t)\,dt}$$

At this point we can show by induction that the solution $\chi_n(\xi)$
grows at most  like $\exp(2A/3\xi^{3/2})\xi^{n/2}$.
So we assume that this holds for $s\le n$ and we show
that it is true for $n+1$. Using the remark and  (\ref{biginduct}) the
induction step is: with $p_n=n/2$
$$\max_{0\le\sigma\le n;\sigma\le t\le n+3}
                \left\{\left[\frac{1}{2}(n+1-t)\right]+
                \frac{(-1)^n+1}{2}+\frac{1}{2}(t-\sigma)+
                p_n\right\}\le p_{n+1}+1$$
which is straightforward.

Finally, we argue that there are two linear independent
formal solutions of this type that can be written in the form
(\ref{expfrm}) which is convenient
for our approach. For this we have to choose
$\chi_0$, which is a solution
of the homogeneous Airy equation such that it does not vanish in $D_{\alpha}$.
Since the Wronskian of the couple Ai($\xi$), Bi($\xi$) is 
a nonzero constant, any combination with real nonzero constants 
of the form 
$C_1$Ai($\xi$)+$i C_2$Bi($\xi$) is an everywhere nonzero 
solution (and the derivative is also nonzero). We can choose
two linear independent solutions in this way, say the ones
for which $\chi_0=$Ai($\xi$)$\pm i$Bi($\xi$). That they
are formally linearly independent with respect to the solutions gotten
from (\ref{ai1}) follows easily, for instance from
the fact that they correspond to different roots of the
characteristic polynomial.
\vskip 0.5cm

To show that there is an actual solution
for each formal solution in this region we first
single out a true solution corresponding to 
the dominant characteristic root and then use it
to reduce the problem to a regular one. Then we show
that they give the expected asymptotic behavior
for the solutions of the original equation.

The ideas are similar to those used in the regular 
case with the exception that extra care is needed
along the degenerate directions.

Choose 
$$\chi_0(\xi)={\mbox{Ai}(\xi)+i\,Bi(\xi)}$$
and consider the  nonzero formal solution that, has the
leading order $\chi_0$.
We proceed as in Step 1, Section 3 to 
construct a rescaled recurrence with respect to the truncation 
of our formal solution. Exactly the same argument as there
shows that the new coefficients $\tilde a$ have smooth asymptotic
series.  

What is new here is that we must provide for the estimate
of the type c) to which end we examine the complete characteristic equation
 $P(\xi,\d;\tilde\lambda)=0$. 
P is a polynomial  in $\tilde\lambda$ (actually it is, to leading order, a
polynomial with constant coefficients) 
$C^{\infty} $ in $\xi$ and $\d$.
P(0,0,$\tilde\lambda$) has a double root $\tilde\lambda=1
$  but by assumption the second derivative does not vanish
so that we can obtain the roots of the polynomial perturbatively.
After series expansion, we obtain:

\be{serpoly}\sum_{j=0}^l\left\{\left[
               a_j(0,0)\left(1+\d\,C_j+\d^2 E_j+O(\d^3)\right)
               +\xi D_x a_j\d^2+O(\d^3)\right]\right.$$$$
\left.          \left(1+j\chi_1(\xi)\d+\left(j\chi_2(\xi)
                \d^2+j(j-1)/2\chi_1^2\right)+O(\d^3)\right)
                \right\}=0  
\ee

\z where we have taken $\tilde\lambda\sim
1+\chi_1(\xi)\d+\chi_2(\xi)\d^2+O(\d^3)$ and

$$C_j=\frac{1}{\chi_0}\left(j\chi'_0+\chi_1\right)$$
$$E_j=\frac{1}{\chi_0^2}\left(j^2\chi_0''\chi_0+
j\chi'_1\chi_0-j\chi_1\chi'_0-\chi_1^2\right)$$

\z Using the relations (\ref{poly}) we get two solutions for
$\chi_1$: $\chi_1=0$   ---actually, as expected,
we get a root $\tilde\lambda_1=1+O(\d^s)$---  and 
$\chi_1=1-2\chi_0'\chi_0^{-1}\delta+O(\d^2)$. We see that in the first order
in $\d$ there is no  root--- crossing, which is not a surprize
since a generic perturbation tends separate  coalescing roots.
The asymptotic series 
are uniformly valid in our domain $D_{\alpha}$. Now
we show that 
\be{smallstep}\left|\frac{\l_2((k+1)\e)-\l_2(k\e)}
              {\l_1(k\e)-\l_2(k\e)}                          \right|\,
                                      <\,\frac{{\mbox{const}}}{|k|+1}
\ee

\z(this explains the condition c) at the begining of this section).

We have 

\z1) 
$$\left|\left(\log\,\chi_0(\xi)\right)'\right|>\ct\sqrt{\xi}+\ct'$$
This is obvious  since it holds asymptotically and the function does not
vanish. Hence
$$\left|\tilde\lambda_2(\d\,k))-\tilde\lambda_1(\d\,k)\right|>\ct\sqrt{\d
|k|+1}$$
2)
$$\left|\left(\log\, h_0(\xi)\right)''\right|<\ct/\sqrt{\xi+1}+\ct'$$

Using the asymptotic series for the $\tilde\lambda_{1,2}$ and the
estimates
 1.) and 2.) we see that

$$\left|\frac{\tilde\lambda_2(\d\,(k+1))-\tilde\lambda_2(\d\,k)}
           {\tilde\lambda_2(\d\,k)-\tilde\lambda_1(\d\,k)}\right|<
            \frac{(\sqrt{\xi+1}+\ct')\d^2}{(\ct\sqrt{\xi}+\ct')\d}<
            \frac{\ct}{|k|+1}$$

The similar estimates for the other roots are better but  this
of course does not improve the overall rate
of convergence.
$\tilde\lambda_1$ =1+O($\d^s$) and can be obviously  made less than 
$1+\ct\,\delta^2$ in  $D_{\alpha}$ so that also $|\tilde\lambda_1((k+1)\d)-
\tilde\lambda_1(k\d)|<1+\ct'\,\delta^2$ which is enough
for our purposes. For $m\ne 1,2$ we can for instance 
use the fact that the derivative of the polynomial at these points
 does not
vanish and settle for  a crude bound
$|\tilde\lambda_m(\d\,(k+1))-\tilde\lambda_m(\d\,k)|<\ct\,\d^2$ which
can be obtained immediately from (\ref{serpoly}).

 Now,  to see that there
is a true solution of the recurrence which is asymptotic to our
formal series starting with Ai+iBi, we only have to repeat 
the same arguments as in the regular case.
\vskip 1cm

The next step is to use this particular solution to lower
the order of the recurrence. We mimic the construction
done in the proof of Lemma \ref{LB} to get a lower 
order recurrence in the variable $d_k$

$$\sum_{s=0}^{l-1}b_s(k\d ,\d )d_{k+s}=0$$
(see (\ref{38})) and
want to check that this new recurrence satisfies the
hypothesis of 
of Lemma \ref{A1}. To leading order, the characteristic
polynomial of the above equation
has no double roots (it now 
has only one root equal to $1$).
\vskip 0.3cm
 
\z As in the regular
case, $b_0(\xi,0)=-\tilde a(\xi,0)=-a_0(0,0)$; 
$b_{l-1}(\xi,0)=-\tilde a_l(\xi,0)=-a_l(0,0)$.
\vskip 0.3cm

\z Noting
that the coefficients of the recurrence (\ref{38}) have asymptotic
series valid throughout the domain (as finite sums
of terms of the series of
$\a_j(\xi,\d)Y_{1,k+j}/Y_{1,k}$), the boundedness
of the coefficients is also trivial.
\vskip 0.3cm
\z As in the proof of Lemma \ref{LB},
the polynomial 
 in $b_j$ has  the same roots as
the polynomial in $\tilde a_j$ (except for the eliminated one) for which we
have already obtained  the estimates
of type c).

Now we have a regular problem for which we know that to each
formal solution there is a genuine solution asymptotic to it.
\vskip 1cm

It remains to check that we can recover the asymptotic
behavior of the solutions of the original recurrence
from those of the reduced one. For the solutions
corresponding to the characteristic roots that
are less than one, exactly the same proof as
in Step 3 of Lemma \ref{LB} works. 
For any formal solution of our original equation
that corresponds to the largest eigenvalue
and which does not vanish, the proof
is the one given in Lemma \ref{A1}.

 In the crossing region however there
might be a special interest in finding a particular solution
which is not of exponential type and
which is small for large
$\xi$ ( the Airy--like solution). To this end,
a slightly different argument is necessary. We can obtain a formal 
solution of the reduced
equation which is, to leading order,

\be{formA}
{{   \frac{\mbox{Bi}(\xi+\d)\mbox{Ai}(\xi)-\mbox{Bi}(\xi)\mbox{Ai}(\xi+\d)}
                {(\mbox{Ai}(\xi)-i\mbox{Bi}(\xi))(\mbox{Ai}(\xi+\d)-i\mbox{Bi}(\xi+\d))}}}
\ee

\z (suggested by computing (\ref{defd}) for two formal
solutions of the original equation, corresponding 
to Ai$\pm$iBi)). The asymptotic 
behavior of (\ref{formA}) for large $\xi$ is 

\be{asyd}
{{\frac{\d}{(\mbox{Ai}(\xi)-i\mbox{Bi}(\xi))^2}(1+\mbox{power series})}}
\ \ (1\ll|\xi|\ll\d^{-1/2})
\ee

Writing the asymptotic  representation of  Ai$\pm$iBi in the 
form 
$$\sim\xi^{-1/4}\exp(2/3\xi^{3/2})(1+series)$$
is sufficient to see this. It is important to
note that there are no  powers of $\xi$ multiplying  the asymptotics
(\ref{asyd});
its leading order does  not vanish and (\ref{asyd}) can 
be written as an exponential of a formal series,
the form required by our arguments.

We now apply the construction in Step 3, Section 3 to
recover the solutions of the initial recurrence.

Using the asymptotic behavior of the Airy functions for large
argument, we get for the reconstructed solution the 
representation

$$Y\sim\xi^{(-1/4)}\exp(\frac{2}{3}\xi^{3/2})
           \sum_{k=j}^{\d^{-\alpha}}(k\d)^{1/2}\exp(-\frac{4}{3}(k\d)^{\frac{3}{2}})
          (1+power\ series)$$
for which the Euler--Maclaurin summation formula gives
the asymptotic representation,

$$\xi^{-1/4}\exp\left(-\frac{2}{3}\xi^{3/2}\right)$$

\z for large positive $\xi$ 

In conclusion there is a true solution of the recurrence
which behaves like the Airy function for positive $k$
(also for negative $k$ when this is properly interpreted)  and 
our argument shows what initial conditions have
to be chosen in order to obtain it. 
For negative $k$ we see that in fact  all the solutions 
corresponding to the largest two eigenvalues are
comparable.

{\bf B.} The exterior region.

Now we want to show that the solutions coming from
the exterior region remain asymptotic to the
true solutions as long as $|x|\gg\e^{2/3}$

The problem that arises here is that
the characteristic polynomial has
a virtually degenerate root for small $x$ and
this leads to a lesser smoothness of the asymptotic
series and ultimately to its collapse at $|x|\sim \e^{2/3}$.
Let $\beta$ be as in Proposition 2.2 and 
define the exterior region by
$$E_{\beta}=\{x:\ |x|>\e^{\beta}\}$$
In what follows we make the following conventions. We write
$$g(\e)=O(\e^{\infty})$$
if $\lim_{\e\rightarrow 0}\e^{-k}\,g(\e)=0$ for all $k$; also
$\cf(x^{\gamma})$ will denote a generic function such that it
together with its derivatives of arbitrary order $s$  satisfy the estimates
\be{have}
\left|\cf(x^{\gamma})^{(s)}(x)\right|<A_{s}+B_{s}x^{\gamma-s}
\ee
with the convention that $|\cf(x^0)|\equiv|\cf(\ln(x))|<A+B|\ln|x||$.

The first step is to study the asymptotic properties
of the formal solutions i.e of the (possibly divergent)
expressions for which
\be{fore}
\sum_{j=0}^l a_j(x;\e)
\exp\left(\sum_{t=0}^{\infty}\Phi_t(x+j\e)\e^t\right)\equiv 0
\ee

We show by induction that $\Phi_t$ and their derivatives behave like
$x^{\frac{3}{2}(1-t)}$ and its derivatives.
We place ourselves in the assumption of genericity of the
crossing which means in particular that ${\frac{\partial
P(\lambda;x)}{\partial\lambda}}>const\sqrt{|x|}$.

\bl{lemme}
i) If in (\ref{fore}) the asymptotic series for the $a_j$ are of the 
form \break $\sum_{k=0}^{\infty}\cf_k(x^{1-\frac{3}{2}k})\e^k$ then, in 
the formal solution (\ref{fore}) we have
$\Phi_t$=$\cf(x^{\frac{3}{2}(1-t)})$. 

ii) The same conclusion is 
true if ${\frac{\partial P(\lambda;x)}{\partial\lambda}}>const>0$
uniformly in $x$ and $a_j\sim \sum_{k=0}^{\infty}\cf_k(x^{\frac{1}{2}-\frac{3}{2}k})\e^k$
\end{lemma}

Note that the coefficients $a_j$ of our initial recurrence
are smoother then it is assumed in i) but this
smoothness does not withstand a rescaling 
as done for  (\ref{37}).

The proof of the Lemma is by induction on $t$.

i) It is easy to see from the eikonal equation that $\Phi_0$=$\cf(x^{3/2})$.
Assume that the conclusions of the lemma hold for all $t'<t$.  After 
a formal series expansion of the exponent of  (\ref{fore}) one gets

\be{xpo}
\hskip 6pt \sum_{j=0}^l a_j(x;\e) \exp\left[\e^t\left(
                                          j\Phi_t^{(s)}+\hskip -20pt
           \sum_{{{\tiny\begin{array}{ccc}s\ge 1\\ t'+s= t+1\end{array}}}}
                \hskip -12pt \frac{j^s}{s!}\Phi_{t'}^{(s)}\right)
                  +\sum_{k=0}^{t}\e^{k-1}\hskip -12pt
            \sum_{{{\tiny\begin{array}{ccc}s\ge 1\\ t'+s\le t+1\end{array}}}} 
                                    \hskip -12pt
                        \frac{j^s}{s!}\Phi_{t'}^{(s)}+O(\e^{t+1})\right]=0
\ee
which, using the induction hypothesis, can be rewritten as

$$
\sum_{j=0}^l\left(\sum\cf_{k;j}(x^{1-\frac{3}{2}k})\e^k\right)
          \exp\left(j\Phi'_0+\e^t\left(j\Phi'_t+\cf(x^{1-\frac{3}{2}t})\right)
          +\right.$$
\be{xpos}
\left.
\sum_{k=1}^{t-1}\e^{k}\cf_k(x^{\frac{1-3k}{2}})+O(\e^{t+1})\right)=0
\ee

After expanding in powers of $\e$ and collecting the term in 
$\e^t$ we obtain the equation for $\Phi_t$ in the form:
\be{eqt}
\sum_{j=0}^l\left(j a_{j;0}(x)\e^{\Phi'_0(x)}\,\Phi_t'(x)+
                    \cf(x^{1-\frac{3t}{2}})\right)=0
\ee
or
\be{eqtt}
 \Phi'_t(x)=
\frac{\cf(x^{1-\frac{3t}{2}})}{\frac{\partial P(\lambda;x)}{\partial
\lambda}}
\ee
where the derivative of the polynomial is evaluated at $\e=0,
{\lambda=\exp(\Phi'_0)}$ thus proving  (i).

For (ii) the same proof works, replacing everywhere
$\cf(x^{1-\frac{3t}{2}})$ with $\cf(x^{\frac{1}{2}-\frac{3t}{2}})$.

\z${\bf\ B1.}$ Rescaling first
the recurrence with respect to the approximate solution we
show that there is a genuine solution corresponding to the maximal
eigenvalue. 
Let $\Psi_m$ be any function 
such that $\Psi_m(x;\e)\sim\sum_0^{\infty}\Phi_{m;t}(x)\e^{t}$
and take
\be{yform}
\tilde Y_{m}:=\exp(\ei\,\Psi_m(x;\e))
\ee

The existence of a solution corresponding to the asymptotics
$\exp(\ei\,\Psi_1(x;\e))$ is again equivalent to a solution
which is $1$ to all orders in $\e$ of the recurrence

\be{28bis}
\sum_{j=0}^l \tilde a_{j,s}(k\epsilon ,\epsilon )C_{k+j}=0
\ee
where $ \tilde{a}_{j,s}(k\epsilon ,\epsilon )={a}_j(k\epsilon
,\epsilon
)\tilde{Y}_{1}((k+j)\epsilon )/
\tilde{Y}_{1}(k\epsilon )$.
The formal solutions of the equation (\ref{28bis}) are
$\tilde Y_{m}/\tilde Y_{1}$. We need the roots of the new
characteristic
polynomial
\be{fullpoly}
\tilde P(\tilde\lambda):=\sum_{j=0}^l \tilde a_{j,s}(k\epsilon ,\epsilon )\tilde\l^j=0
\ee

\z It is easy to see that the polynomial (\ref{fullpoly}) has a root 
which is $1$ to all orders in $\e$. Let now $G(x;\e)$
be one of the differences $G_m(x;\e)=\Psi_m(x;\e)-\Psi_1(x;\e)$. We have,
$$\sum_{j=0}^l \tilde a_{j,s}(k\epsilon ,\epsilon )\exp(\ei G(x+j\e;\e))=
o(\e^t)$$

\z for all $t$ and so, after series expansion

\be{polyexp}
\sum_{j=0}^l \tilde a_{j,s}(k\epsilon ,\epsilon )
              \lambda_0^j (1+\e H_j(x;\e))=o(\e^{\infty})
\ee

\z where $H_j(x;\e)$ are some smooth functions of $x,\e$ and
$\lambda_0:=\exp(G_x(x;\e))$. Using Lemma \ref{lemme}
and the genericity assumptions it is not
difficult  to see that

\be{estH}
|H_j(x;\e)|<\ct\,|x|^{-1/2}
\ee

If we look for solutions
of the characteristic polynomial (\ref{fullpoly}) in the form $\lambda_0+\gamma$
we get

$$0=\tilde P(\tilde\lambda)=\sum_{j=0}^{l} P^{(j)}(\lambda_0)
\gamma^{j}$$

\z (where the derivatives are taken with respect to $\lambda$). 
Using (\ref{polyexp}) we obtain  $\gamma$  as the
unique small solution of the equation

\be{coe}
\gamma=\frac{\e\sum_{j=0}^{l}\tilde a_j(x;\e)\lambda_0^j H_j(x;\e)}
              {\sum_{j=1}^{l} P^{(j)}(\lambda_0)\gamma^{j-1}}
\ee

which is a contraction for small enough $\e$ (and small $\gamma$ ) 
in the 
region $|x|>\e^{\beta}$ as it
is easy to check. We then obtain from (\ref{estH}) and (\ref{coe})

\be{estga}
|\gamma|<\ct\frac{\e}{|x|}
\ee
again valid for $|x|>\e^{\beta}$.

We shall also need estimates for $\gamma(x+\e)-\gamma(x)$.
Differentiating  (\ref{coe}) with respect to $x$ and using 
Lemma \ref{lemme}

\be{estdiff}
|\gamma(x+\e)-\gamma(x)|<\ct\,\e\,|x|^{-2}
\ee

Now we proceed as in the regular case in rewriting
the recurrence in matrix form and evaluating the terms
in the product (\ref{32}). In the  matrices 
$T:=\tilde \Lambda_{k}^{-1}\tilde \Lambda_{k-1}$
the off--diagonal elements are estimated by $T_{mn}<\ct\,\e|x|^{-1}$
and for the diagonal elements we have $T_{mm}=1+O(\e/x)$. Indeed,

$$T_{mm}=\prod_{n\ne m}\left(1+\frac{\l_m((p-1)\e,\e)-{\l_n(p\e,\e)}}
{\l_m(p\e,\e)-\l_n(p\e,\e)}\right)$$
Since by the assumption of genericity
the roots of the polynomial are separated 
by at least $\ct\sqrt{(|x|)}$ each  term in the product above can be estimated by

$$1+\ct\,\e(\l_0'(x)+\gamma'(x))|x|^{-1/2}+O\left(\e(\l_0'(x)+\gamma'(x))|x|^{-1/2}\right)
$$$$<1+\ct\left(\,\e|x|^{-1}+\e^2|x|^{-5/2}\right)<1+\ct\,\e|x|^{-1}$$

\z in our
region $E_{\beta}$. The nondiagonal terms are estimated in a very 
similar way.

We derive the estimate 
 $\|T(x)-{\bf I}\|<K\e/x$ for some constant $K$. Assume for 
definiteness that we are on the left of the crossing point. We get,
$$\left\|\prod_{k; -k\e>\e^{\beta}}T(k\e)\right\|<\ct\,\e^{-K/3}$$

Finally we have to control  the product
of the norms of the diagonal matrices
$D_{k}$. Since they all have one   eigenvalue 
equal to $1$ to all orders
in $\e$ and for $i>2$ $|\tilde\lambda_i(x;\e)|<1-\ct$ the only
nontrivial contribution comes from $\tilde\lambda_2$ and this
only if $\tilde\lambda_1$ and $\tilde\lambda_2$ have the
same modulus to leading order in $\e$.
Referring
to the decomposition $\tilde\lambda=\l_0+\gamma$ we have in this
case, using Euler-Maclaurin summation formula,
$$\prod_{k; -k\e>\e^{\beta}}|\l_0(k\e)|=\exp\left(\ei\sum_k\Re((\Psi_2)_x(k\e;\e)-(\Psi_1)_x(k\e;\e))
\right)\sim$$ 
$$\exp\left\{\Re\left(\Phi_{2;1}(-\e^{2/3})-\Phi_{1;1}(-\e^{2/3})\right)\right\}<\e^{-const}$$
Also,
$$\prod_{k; -k\e>\e^{\beta}}\left|1+\frac{\gamma(k\e)}{\l_0(k\e)}\right|
<\prod_{k;
-k\e>\e^{\beta}}\left|1+\frac{\ct}{|k|}\right|<\e^{-const}$$

\z so that also  $\prod\tilde\l$ is less than $\e^{-const}$. 
At this point the same arguments as in the regular case show 
that there is a true solution behaving asymptotically as
the formal solution corresponding to the largest eigenvalue.

${\bf B2}$

Let $Y_1$ be a solution of the recurrence relation 
such that $Y_1\sim\exp\left(\ei\sum\Phi_{1;t}\e^t\right)$ in
$E_{\beta}$. We now follow the same steps that led to
equations (\ref{37}) and (\ref{38}).

It is a matter of straightforward induction to derive from
Lemma \ref{lemme} that the 
coefficients $\tilde a_j$ have the behavior

$$\tilde a_j(x;\e)\sim a_{j;0}(x)e^{j\,\Phi'_0(x)}+
                     \sum_{k\ge
1}\cf_{j;k}(x^{\frac{1}{2}-\frac{3}{2}k})\e^k$$

\z and then clearly 
\be{estimb}
b_j(x;\e)\sim\sum_{k\ge
0}\cf_{j;k}(x^{\frac{1}{2}-\frac{3}{2}k})\e^k
\ee
It is also easy to check that the recurrence (\ref{38})
is now nondegenerate in the sense that:
\be{nondg}
\inf_{E_{\beta}}\left\{|b_0(k\e)|,|b_{l-1}(k\e)|,\frac{1}{|b_j(k\e)|}\right\}>{\mbox{const}}>0
\ee

\z 
and the characteristic polynomial of the new recurrence does
not have coalescing  roots (the root $\l=1$ of
(\ref{37}) has been eliminated in the reduction):
\be{separmod} \inf_{E_{\beta}}\{\left|\,|\lambda_m(x)|
-|\lambda_n(x)|\,\right|\}>{\mbox{const}}>0\ \ ({\mbox{for }}m\ne n)
\ee

We are now left with a problem of the following type.
Taking a recurrence of the form

\be{nonsmoothrec}
\sum_{j=0}^l a_j(x;\e)y_{k+j}=0
\ee
 
under the following conditions:
\be{as1}
\tilde a_j(x;\e)\sim 
                     \sum_{k\ge 0}\cf_{j;k}(x^{\frac{1}{2}-\frac{3}{2}k})\e^k
\ee

\be{nondgpl}\inf_{E_{\beta}}\left\{|a_0(k\e)|,|a_{l-1}(k\e)|,\frac{1}{|a_j(k\e)|}\right\}>{\mbox{const}}>0
\ee
\be{separmodl} \inf_{E_{\beta}}\{\left|\,|\lambda_m(x)|
-|\lambda_n(x)|\,\right|>{\mbox{const}}>0\ (m\ne n)\}
\ee

where $\lambda_m(x)$ are the roots of the polynomial

\be{redpoly}
P(\lambda):=\sum_{j=0}^l a_{j,0}(k\epsilon )\l^j=0
\ee
we want to show that
\bl{almostreg} Given a formal solution to (\ref{nonsmoothrec}):

$$\exp\left(\ei\sum_{t=0}^{\infty}\e^t\Phi_t(x)\right)$$

\z where $\Phi_t(x)=\cf_t(x^{\frac{3}{2}(1-t)})$,
 there is a solution asymptotic to it 
for $0<x\in E_{\beta}$ (and correspondingly one 
when $x$ negative).
\end{lemma}

Proof: induction on $l$.

a) We show that we can find a solution corresponding to the root
the has the largest modulus (this will simultaneously
prove the lemma for $l=1$). All the arguments in B1
above apply here. Actually, now  we could get
some better estimates since we do not have small
denominators in  (\ref{estga}), (\ref{estdiff}) 
and in the estimates of the matrices but this would not
affect the final result.

b) We assume that the conclusion is true for all recurrences of 
order less than $l-1$ and show it holds for recurrences
of order $l$. By the arguments above, there is
a true solution asymptotic to the formal solution defined
by the maximum eigenvalue.
Using it to reduce the order
of the recurrence we obtain an order--$l-1$ scheme, 
which satisfies the conditions of Lemma \ref{almostreg} as
it is easy to check 
and for which we thus know the asymptotic behavior of the
solutions. It remains to verify that they can be
used to produce  solutions of the higher--order recurrence
with the stated asymptotic behavior. For definiteness
we study the subregion $x<\ -\e^{\beta}$. All
the arguments in Step 3, Section 3 apply if we
take $k_1$ to be $-\e^{\beta}$. The only 
change is that in (\ref{endest}) $\Phi_{m;1}(x)$
are not uniformly bounded. Instead, using
Lemma \ref{lemme} we get
$|\Phi_{j;1}|=\cf(\ln(|x|))<K|\ln\,\e|$ for
some fixed constant $K$ so the RHS of (\ref{endest})
changes to $Y_{m;k}\left(1+O(\e^{s-2-K})\right)$

The conclusion is that after the first reduction we end
up with a recurrence that is nondegenerate in the sense
of Lemma \ref{almostreg} and for which we can control the
small--$\e$ behavior of the solutions. Now, 
the reconstruction of the solutions of the original recurrence
from the solutions of the reduced one amounts
to merely repeating without any significant 
change the construction and estimates in 
part b) above. At this point in the proof it is clear
that if the crossing roots are not the largest, on can 
reduce the order of the recurrence to the actual level at which the roots
cross, and then apply the arguments above.

\end{section}
$$ $$
$$ $$

\eject
\bibliography{siam}
\bibliographystyle{siam}
\par [1] Percy Deift, Ken McLaughlin, A Continuum Limit of the
Toda Lattice, in preparation
\par [2] Courant, R.,  Hilbert, D{\it Methods of Mathematical Physics}
John Wiley, 1962
\par [3] Erdelyi, A. {\it Higher transcendental functions} McGraw
Hill,1953
\par [4] Borel, E, {\it Le{\c c}ons sur les s\'eries div{\'e}rgentes},
  Gauthier-Villars
(1928)
\par [5] Remmert, R. {\it Theory of complex functions}, 
Springer-Verlag, 1991
\end{document}